\documentstyle[12pt]{article}
\def\a{\alpha}      \def\b{\beta}          \def\g{\gamma}
\def\d{\delta}              
      \def\th{\theta}

\def\CC{{\bf C}\!\!\!{\rm l}}

\def\ad{{\rm ad}}
\def\wt{\widetilde}
\def\gu{\widetilde{g[u]}}
\newcommand{\bn}{\begin{equation}}
\newcommand{\ed}{\end{equation}}
\newtheorem{definition}{Definition}[section]
\newtheorem{proposition}{Proposition}[section]
\newtheorem{theorem}{Theorem}[section]

\normalbaselines
\begin{document}
\begin{center}
\vspace*{1.0cm}

{\LARGE{\bf Drinfeldians}}

\vskip 1.5cm

{\large {\bf Valeriy N. Tolstoy}}

\vskip 0.5 cm

Institute of Nuclear Physics\\
Moscow State University\\
119899 Moscow\&Russia\footnote
{\it e-mail: tolstoy@anna19.npi.msu.su}

\end{center}

\vspace{1 cm}

\begin{abstract}
We construct two-parameter deformation of an universal
enveloping algebra $U(g[u])$ of a polynomial loop algebra $g[u]$, 
where $g$ is a finite-dimensional complex simple Lie algebra 
(or superalgebra). This new quantum Hopf algebra called 
the Drinfeldian $D_{q\eta}(g)$ can be considered as 
a quantization of $U(g[u])$ in the direction of a classical 
r-matrix which is a sum of the simple rational and trigonometric
r-matrices. The  Drinfeldian $D_{q\eta}(g)$ contains $U_{q}(g)$ 
as a Hopf subalgebra, moreover $U_{q}(g[u])$ and $Y_{\eta}(g)$ are
its limit quantum algebras when the $D_{q\eta}(g)$ deformation 
parameters $\eta$ goes to $0$ and $q$ goes to $1$,
respectively. These results are easy generalized to a supercase, 
i.e. when $g$ is a finite-dimensional contragredient simple 
superalgebra.
\end{abstract}

\vspace{1 cm}

\section{Introduction}
As it is well known, an universal enveloping algebra
$U(g[u])$ of a polynomial loop (current) Lie algebra $g[u]$ ,
where $g$ is a finite-dimensional complex simple  Lie algebra,
admits two type deformations: a trigonometric deformation
$U_q(g[u])$ and a rational deformation or
Yangian $Y_{\eta}(g)$ \cite{D1,D2}.
(In the case $g=sl_n$ there also exists an elliptic quantum
deformation of $U(sl_n[u])$).
The algebras $U_q(g[u])$, and $Y_{\eta}(g)$ are quantizations
of $U(g[u])$ in the direction of the simplest trigonometric
and rational solutions of the classical Yang-Baxter equation
over $g$, respectively. These deformations are one-parameter ones.
It turns out that $U(g[u])$ also admits two-parameter
deformation which is called the Drinfeldian or
the rational-trigonometric quantum algebra.
The Drinfeldian $D_{q\eta}(g)$ is a quantization
of $U(g[u])$ in the direction of a classical r-matrix which
is a sum of the simplest rational and trigonometric r-matrices.
The  Drinfeldian $D_{q\eta}(g)$ contains $U_{q}(g)$ as a Hopf
subalgebra, and  $U_{q}(g[u])$ and $Y_{\eta}(g)$ are its limit
quantum algebras when the deformation parameters of $D_{q\eta}(g)$
$\eta$ goes to $0$ and $q$ goes to $1$, respectively.
These results are easy generalized to a supercase, i.e. when
$g$ is a finite-dimensional contragredient simple superalgebra.

\setcounter{equation}{0}
\section{Quantum algebra $U_q(\gu)$}

Let $g$ be a finite-dimensional complex simple Lie algebra
of a rank $r$ with a standard Cartan matrix
$A=(a_{ij})_{i,j=1}^r$, with a system of simple roots
$\Pi:= \{\alpha_1,\ldots, a_r\}$ and a maximal positive
root $\theta$, and with a Chevalley basis
\{$h_{\alpha_i}^{}$, $e_{\pm\alpha_i}'$, $(i=1,2,\ldots, r)$\}.
Let $g[u]$ be a polynomial loop algebra (or a Lie algebra
of polynomial currents over $g$), and $\widetilde{g[u]}$
be a "central extension" of $g[u]$:
$\widetilde{g[u]}\simeq g[u]\oplus\CC\,c$, where $c$
is a central element.\footnote{More correctly, the element $c$
is a central element of a central extension of the total loop
algebra $g[u,u^{-1}]$.}
The Lie algebra $\widetilde{g[u]}$
(and its universal enveloping algebra $U(\widetilde{g[u]})$)
is generated by the Chevalley basis of $g$ and the affine
element $e_{(\delta-\theta)}':=ue_{-\theta}^{}$ and
$h_{\delta}^{}:=c$ with the following defining relations:
\begin{eqnarray}
&&[h_{\delta}^{},{\rm everything}]=0~,\qquad\qquad\qquad
[h_{\alpha_{i}}^{}, h_{\alpha_{j}}^{}]=0~,
\label{A1}
\\[7pt]
&&[h_{\alpha_{i}}^{},e_{\pm\alpha_{j}}']=
\pm(\alpha_{i},\alpha_{j})e_{\pm\alpha_{j}}'~,\qquad\quad
[e_{\alpha_{i}}',e_{-\alpha_{j}}']=
\delta_{ij}h_{\alpha_{i}}^{}~,
\label{A2}
\\[7pt]
&&(\ad\,e_{\pm\alpha_{i}}')^{1-a_{ij}}e_{\pm\alpha_{j}}'=0
\qquad\qquad\qquad{\rm for}\,\,i\neq j~,
\label{A3}
\\[7pt]
&&[h_{\alpha_i}^{},e_{\delta-\theta}']=-
(\alpha_i,\theta)\,e_{\delta-\theta}',\qquad\qquad
[e_{-\alpha_i}',e_{\delta-\theta}']=0~,
\label{A4}
\\[7pt]
&&(\ad\,e_{\alpha_i}')^{1-a_{i0}}e_{\delta-\theta}'=0,
\qquad\qquad\qquad
\big(a_{i0}=-2(\alpha_i,\theta)/(\alpha_i,\alpha_i)\big)~,
\label{A5}
\\[7pt]
&&[[e_{\alpha_i}',e_{\delta-\theta}'],e_{\delta-\theta}']=0
\qquad\qquad\qquad\;\;{\rm for}\; g\ne sl_2 \;\;
{\rm and}\quad(\alpha_i,\theta)\ne 0~,\qquad
\label{A6}
\\[7pt]
&&[[[e_{\alpha}',e_{\delta-\alpha}'],
e_{\delta-\alpha}'],e_{\delta-\alpha}']=0
\qquad\qquad{\rm for}\; g=sl_2~.
\label{A7}
\end{eqnarray}
Here $"\ad"$ is the adjoint action of $\widetilde{g[u]}$ in
$\wt{g[u]}$, i.e. $(\ad\,x)y=[x,y]$ for any $x,y\in \gu$.
The relations (\ref{A6}) relate to the case $g\ne sl_2$,
and the relation (\ref{A7}) belongs to the case $g=sl_2$
(in this case $\theta=\alpha_{r=1}$ and we set
$\alpha:=\alpha_{r=1}$).

{\it Remark}. The defining relations for $\gu$
can be obtained from  defining relations of the
corresponding non-twisted affine Lie algebra $\hat{g}$
by removing relations with a negative affine root vector
$e_{-\delta+\theta}'$.

Let $U_{q}(g)$ be a standard q-deformation of
the universal enveloping algebra $U(g)$ with
Chevalley generators $k_{\alpha_i}^{\pm 1}$, $e_{\pm\alpha_i}^{}$
$(i=1,2,\ldots, r)$ and with the defining relations
\vspace{-3pt}
\begin{eqnarray}
&&[k_{\alpha_i}^{},k_{\alpha_j}^{}]=0~,\qquad\qquad\qquad
k_{\alpha_i}^{}e_{\pm\alpha_j}^{}k^{-1}_{\alpha_i}=
q^{\pm(\alpha_i,\alpha_j)}e_{\pm\alpha_j}^{}~,
\label{A8}
\\[7pt]
&&[e_{\alpha_i}^{},e_{-\alpha_i}^{}]=
\frac{k_{\alpha_i}^{}-k_{\alpha_i}^{-1}}{q-q^{-1}}~,\qquad
({\rm ad}_{q}e_{\pm\alpha_{i}}^{})^{1-a_{ij}}e_{\pm\alpha_{j}^{}}=0
\quad\;{\rm for}\,\, i\neq j~,\qquad
\label{A9}
\end{eqnarray}
where $(\ad_qe_{\beta}^{})e_{\gamma}^{}$ is the q-commutator:
\begin{equation}
(\ad_qe_{\beta}^{})e_{\gamma}^{}:= [e_{\beta}^{},e_{\gamma}^{}]_q:=
e_{\beta}^{} e_{\gamma}^{}-q^{(\beta,\gamma)}e_{\gamma}^{}e_{\beta}^{}\ .
\label{A10}
\end{equation}
A Hopf structure of $U_{q}(g)$ is given the following formulas
for a comultiplication $\Delta_{q}$, an antipode $S_{q}$,
and a co-unite $\varepsilon_{q}$:
\begin{eqnarray}
&&\Delta_{q}(k_{\alpha_i}^{\pm 1})=
k_{\alpha_i}^{\pm 1}\otimes k_{\alpha_i}^{\pm 1}~,
\qquad\quad\Delta_{q}(e_{\alpha_i}^{})=
e_{\alpha_i}^{}\otimes 1+k_{\alpha_i}^{-1}\otimes e_{\alpha_i}^{}~,
\nonumber
\\[4pt]
&&\Delta_{q}(e_{-\alpha_i}^{})=
e_{-\alpha_i}^{}\otimes k_{\alpha_i}^{}+1\otimes e_{-\alpha_i}^{}~,
\label{A11}
\\[7pt]
&&S_{q}(k_{\alpha_i}^{\pm 1})=k_{\alpha_i}^{\mp 1}~,\qquad
S_{q}(e_{\alpha_i}^{})=-k_{\alpha_i}^{}e_{\alpha_i}^{}~,\qquad
S_{q}(e_{-\alpha_i}^{})=-e_{-\alpha_i}^{}k_{\alpha_i}^{-1}~,
\qquad\quad
\label{A12}
\\[7pt]
&&\varepsilon_{q}(k_{\alpha_i}^{\pm 1})=1~,\qquad\quad
\varepsilon_{q}(e_{\pm\alpha_i}^{})=0~.
\label{A13}
\end{eqnarray}
\begin{definition}
The quantum algebra $U_q(\gu)$ (or a q-deformation of
$U(\gu)$) is generated (as an associative algebra) by the algebra
$U_{q}(g)$ and the elements $e_{\delta-\theta}^{}$,
$k_{\delta}^{\pm 1}$ with the relations:
\begin{eqnarray}
&&[k_{\delta}^{\pm 1},{\rm everything}]=0~,\qquad\quad
k_{\alpha_i}^{}e_{\delta-\theta}^{}k^{-1}_{\alpha_i}=
q^{-(\alpha_i,\theta)}e_{\delta-\theta}^{}~,
\label{A14}
\\[7pt]
&&[e_{-\alpha_i}^{},e_{\delta-\theta}^{}]=0~,\qquad\qquad\quad\;\;
(\ad_{q}e_{\alpha_i}^{})^{n_{i0}}e_{\delta-\theta}^{}=0~,
\label{A15}
\end{eqnarray}
where $n_{i0}=1+2(\a_{i},\theta)/(\a_{i},\a_{i})$, and
\begin{eqnarray}
&&[[e_{\alpha_i}^{},e_{\delta-\theta}^{}]_{q},
e_{\delta-\theta}^{}]_{q}=0
\qquad\qquad\quad{\rm for}\; g\ne sl_2\;
{\rm and}\;(\alpha_i,\theta)\ne 0~,\qquad
\label{A16}
\\[7pt]
&&[[[e_{\alpha}^{},e_{\delta-\alpha}^{}]_{q},
e_{\delta-\alpha}^{}]_{q},e_{\delta-\alpha}^{}]_{q}=0
\qquad{\rm for}\,\, g=sl_2~.
\label{A17}
\end{eqnarray}
The Hopf structure of $U_{q}(\gu)$ is defined by
the formulas $\Delta_{q}(x)=\Delta_{q}(x)$,
$S_{q}(x)=S_{q}(x)$ ($x\in U_{q}(g)$), and
$\Delta_{q}(k_{\delta}^{\pm})=
k_{\delta}^{\pm1}\otimes k_{\delta}^{\pm1}$,
$S_{q}(k_{\delta}^{\pm1})=k_{\delta}^{\mp1}$.
The comultiplication, the antipode and the co-unite of 
the element $e_{\delta-\alpha}^{}$  are given by
\begin{eqnarray}
&&\Delta_{q}(e_{\delta-\theta}^{})=e_{\delta-\theta}^{}\otimes 1+
k_{\delta-\theta}^{-1}\otimes e_{\delta-\theta}^{}~,
\label{A18}
\\[7pt]
&&S_{q}(e_{\delta-\theta}^{})=
-k_{\delta-\theta}^{}e_{\delta-\theta}^{}~,\qquad\qquad
\varepsilon(e_{\delta-\theta})=0~.
\label{A19}
\end{eqnarray}
Here we put
\begin{equation}
k_{\delta-\theta}^{}=k_{\delta}^{}k_{\alpha_1}^{-n_1}
k_{\alpha_2}^{-n_2}\cdots k_{\alpha_r}^{-n_r}
\label{A20}
\end{equation}
if $\theta=n_1^{}\alpha_1^{}+n_2^{}\alpha_2^{}+\cdots \,+
n_r^{}\alpha_r^{}$.
\end{definition}
{\it Remark}. The defining relations for $U_{q}(\gu)$
can be obtained from  defining relations of the
corresponding quantum non-twisted affine algebra $U_{q}(\hat{g})$
by removing relations with the negative affine root vector
$e_{-\delta+\theta}^{}$.

It is easy to check the following result.
\begin{proposition}
There is a one-parameter group of Hopf algebra automorphisms
${\cal T}_a$ of $U_q(\gu)$, $a\in \CC$, given by
\begin{eqnarray}
&&{\cal T}_a(k_{\d}^{\pm 1})=k_{\d}^{\pm 1}~,\qquad
{\cal T}_a(k_{\a_{i}}^{\pm 1})=k_{\a_{i}}^{\pm 1}~,
\nonumber
\\[4pt]
&&{\cal T}_a (e_{\pm\alpha_i}^{})=e_{\pm\alpha_i}^{}~,\qquad
{\cal T}_a(e_{\delta-\theta}^{})=a\,e_{\delta-\theta}^{}~.
\label{A21}
\end{eqnarray}
\label{PA1}
\end{proposition}

\vspace{-15pt}
\setcounter{equation}{0}
\section{Drinfeldian $D_{q\eta}(g)$}

Here we keep the notations of the previous Section and begin
with the following important definition.
\begin{definition}
The Drinfeldian $D_{q\eta}(g)$ is generated $($as an associative
algebra over $C\!\!\!\!I\,[[\eta]]$$)$ by the algebra $U_{q}(g)$ and
the elements $\xi_{\delta-\theta}^{}$, $k_{\delta}^{\pm 1}$ with
the relations:
\vspace{-2pt}
\begin{equation}
[k_{\delta}^{\pm 1},{\rm everything}]=0~\qquad
k_{\alpha_i}^{}\xi_{\delta-\theta}^{}k^{-1}_{\alpha_i}=
q^{-(\alpha_i,\theta)}\xi_{\delta-\theta}^{}~,
\label{D1}
\end{equation}
\begin{equation}
[e_{-\alpha_i}^{},\xi_{\delta-\theta}^{}]=
\tau\,[e_{-\alpha_i}^{},\tilde{e}_{-\theta}^{}],\qquad
({\rm ad}_{q}e_{\alpha_i}^{})^{n_{i0}^{}}\xi_{\delta-\theta}^{}=
\tau\,({\rm ad}_{q}e_{\alpha_i}^{})^{n_{i0}^{}}
\tilde{e}_{-\theta}^{}
\label{D2}
\end{equation}
for $n_{i0}^{}=1+2(\alpha_i,\theta)/(\alpha_i,\alpha_i)$~, and
\begin{eqnarray}
[[e_{\alpha_i}^{},\xi_{\delta-\theta}^{}]_{q},
\xi_{\delta-\theta}^{}]_{q}
&\!\!\!\!=&\!\!\!\!\!
-\tau^{2}[[e_{\alpha_i}^{},\tilde{e}_{-\theta}^{}]_{q},
\tilde{e}_{-\theta}^{}]_{q}
\nonumber
\\
&&\!\!\!\!
+\tau\,[[e_{\alpha_i}^{},\tilde{e}_{-\theta}^{}]_{q},
\xi_{\delta-\theta}]_{q}+
\tau\,[[e_{\alpha_i}^{},\xi_{\delta-\theta}^{}]_{q},
\tilde{e}_{-\theta}^{}]_{q}
\end{eqnarray}
\label{D3}
for $g\ne sl_2$ and $(\alpha_i,\theta)\ne 0$,
\begin{eqnarray}
&&\!\!\!\!\!\!\!\!\!\!\!\!\!\!\!\!\!\!
[[[e_{\alpha}^{},\xi_{\delta-\alpha}^{}]_{q},
\xi_{\delta-\alpha}^{}]_{q},\xi_{\delta-\alpha}^{}]_{q}=\tau^{3}
[[[e_{\alpha}^{},\tilde{e}_{-\alpha}^{}]_{q},
\tilde{e}_{-\alpha}^{}]_{q},\tilde{e}_{-\alpha}^{}]_{q}\qquad
\nonumber
\\
&&-\tau^{2}[[[e_{\alpha}^{},\tilde{e}_{-\alpha}^{}]_{q},
\tilde{e}_{-\alpha}^{}]_{q},
\xi_{\delta-\alpha}^{}]_{q}-\tau^{2}[[[e_{\alpha}^{},
\tilde{e}_{-\alpha}^{}]_{q},
\xi_{\delta-\alpha}^{}]_{q},\tilde{e}_{-\alpha}^{}]_{q}
\nonumber
\\
&&-\tau^{2}[[[e_{\alpha}^{},\xi_{\delta-\alpha}^{}]_{q},
\tilde{e}_{-\alpha}^{}]_{q},\tilde{e}_{-\alpha}^{}]_{q}
+\tau\,[[[e_{\alpha}^{},\tilde{e}_{-\alpha}^{}]_{q},
\xi_{\delta-\alpha}^{}]_{q},\xi_{\delta-\alpha}^{}]_{q}
\nonumber
\\
&&+\tau\,[[[e_{\alpha}^{},\xi_{\delta-\alpha}^{}]_{q},
\tilde{e}_{-\alpha}^{}]_{q},\xi_{\delta-\alpha}^{}]_{q}
+\tau\,[[[e_{\alpha}^{},\xi_{\delta-\alpha}^{}]_{q},
\xi_{\delta-\alpha}^{}]_{q},\tilde{e}_{-\alpha}^{}]_{q}
\label{D4}
\end{eqnarray}
for $g=sl_2$. The Hopf structure of $D_{q\eta}(g)$ is defined by
the formulas $\Delta_{q\eta}(x)=\Delta_{q}(x)$,
$S_{q\eta}(x)=S_{q}(x)$ $(x\in U_{q}(g))$ and 
$\Delta_{q}(k_{\delta}^{\pm})=
k_{\delta}^{\pm1}\otimes k_{\delta}^{\pm1}$,
$S_{q}(k_{\delta}^{\pm1})=k_{\delta}^{\mp1}$.
The comultiplication and the antipode of $\xi_{\delta-\alpha}^{}$
are given by
\begin{eqnarray}
\Delta_{q\eta}(\xi_{\delta-\theta}^{})\!\!\!&=&\!\!\!
\xi_{\delta-\theta}^{}\otimes 1+k_{\delta-\theta}^{-1}
\otimes \xi_{\delta-\theta}
\nonumber
\\
&&\!\!\!+ a\left(\Delta_{q}(\tilde{e}_{-\theta}^{})
-\tilde{e}_{-\theta}^{}\otimes 1-
k_{\delta-\theta}^{-1}\otimes\tilde{e}_{-\theta}^{}\right),
\label{D5}
\end{eqnarray}
\begin{equation}
S_{q\eta}(\xi_{\delta-\theta}^{})=
-k_{\delta-\theta}^{}\xi_{\delta-\theta}^{}
+a\left(S_{q}(\tilde{e}_{-\theta}^{})+
k_{\delta-\theta}^{}\tilde{e}_{-\theta}^{}\right)~.
\label{D6}
\end{equation}
Here $\tau:=\eta/(q-q^{-1})$,
$({\rm ad}_{q}e_{\beta}^{})e_{\gamma}^{}=
[e_{\beta}^{},e_{\gamma}^{}]_{q}$,
and the vector  $\tilde{e}_{-\theta}^{}$ is any $U_{q}(g)$
element of the  weight $-\theta$, such that
$g\ni\lim_{q\to1}\tilde{e}_{-\theta}^{}\neq 0$.
\end{definition}
The right-hand sides of the relations (\ref{D2})-(\ref{D6}) are
nonsingular at $q=1$.
\begin{theorem}
(i) The Drinfeldian $D_{q\eta}(g)$ is a two-parameter
quantization of $U(\widetilde{g[u]})$ in the direction of
a classical r-matrix which is a sum of the simplest rational
and trigonometric r-matrices.
\\
\noindent
(ii) The Hopf algebra $D_{q=1,\eta}(g)$
is isomorphic to the Yangian $Y_{\eta}'(g)$
(with the additional central element $c=h_{\d}$).
Moreover, $D_{q\eta=0}(g)=U_q(\widetilde{g[u]})$.
\label{DT1}
\end{theorem}

{\it Remark}. Since the defining relations for $D_{q\eta}(g)$
and $U(\widetilde{g[u]})$ in terms of the Chevalley basis differ 
only in the right-hand sides of the relations 
(\ref{D2})-(\ref{D4}), therefore the Dynkin diagram of $g[u]$ 
can be also used for classification of the Drinfeldian 
$D_{q\eta}(g)$ and the Yangian $Y_{\eta}(g)$.

An analog of Proposition \ref{PA1} is the following result.
\begin{proposition}
There is a one-parameter group of Hopf algebra automorphisms
${\cal T}_a$ of $D_{q\eta}(g)$, $a\in {\CC}$\,, given by
\begin{eqnarray}
&&
{\cal T}_{a}(k_{\d}^{\pm1})=k_{\d}^{\pm1}~,\qquad
{\cal T}_{a}(k_{\a_{i}}^{\pm1})=k_{\a_{i}}^{\pm1}~,\qquad
{\cal T}_{a}(e_{\pm\a_{i}}^{})=e_{\pm\a_{i}}^{}~,\qquad\qquad
\nonumber
\\[5pt]
&&
{\cal T}_{a}(\xi_{\d-\th}^{})=
\left(1-(q-q^{-1})a\right)\xi_{\d-\th}^{}+
\eta \,a\,\tilde{e}_{-\th}^{}~.
\label{D7}
\end{eqnarray}
\label{PD1}
\end{proposition}

\vspace{-10pt}
In the next section we give an explicit description of
the right-hand of the relations (\ref{D2})-(\ref{D6})
for the Drinfeldians $D_{q\eta}(g)$ and the Yangians
$Y_{\eta}(g)$ of the Lie algebras $g$ of rank 2.

\setcounter{equation}{0}
\section{Drinfeldians and Yangians over Lie algebras of rank 2}

Explicit description of the Drinfeldians $D_{q\eta}(g)$ and
the Yangians $Y_{\eta}(g)$  for  the cases $g=sl_2\;and\;sl_3$
were given in \cite{T1}-\cite{T3}. Here we consider the cases
$g=C_{2}(\simeq B_2)\;and\;G_{2}$.

\vspace{5pt}
\noindent
{\large {\it 1. The Drinfeldian $D_{q\eta}(C_{2})$ and
the Yangian $Y_{\eta}(C_{2})$}}.
\noindent
In the case of the Lie algebra $g=C_{2}$ there are two
positive simple roots $\a$ and $\b$, and the maximal
positive root is $\theta=\a+2\b$.

As we already noted the Drinfeldian $D_{q\eta}(g)$ of
the Yangian $Y_{\eta}(g)$ can be characterized the Dynkin diagram
of the corresponding non-twisted Kac-Moody affine Lie algebra
$g^{(1)}$. In the case $g=C_{2}$ the Dynkin diagram  of the
corresponding affine Lie algebra $C_{2}^{(1)}$ is presented by
the picture \cite{K}

{\thicklines
\vskip 15pt
\begin{picture}(500,20)
\put(163,9){\line(1,0){31}}
\put(206,9){\line(1,0){31}}
\put(160,6){\circle{8}}
\put(200,6){\circle{8}}
\put(240,6){\circle{8}}
\put(192,13){\line(1,-2){5}}
\put(192,-1){\line(1,2){5}}
\put(203,9){\line(1,-2){5}}
\put(203,3){\line(1,2){5}}
\put(163,3){\line(1,0){31}}
\put(206,3){\line(1,0){31}}
\put(150,14){$\d\!-\!\th$}
\put(198,14){$\b$}
\put(238,14){$\alpha$}
\end{picture}}
\vskip 8pt
\centerline{\footnotesize Fig.1. Dynkin diagram of the Lie
algebra $C_{2}^{(1)}$}

\vskip 10pt
\noindent
The quantum Hopf algebra $U_q(C_{2})$ is generated by
the elements $k_{\a}^{\pm 1}$, $k_{\b}^{\pm 1}$, $e_{\pm\a}^{}$,
$e_{\pm\b}^{}$ with the defining relations (\ref{A8})-(\ref{A13})
(see also details in \cite{KT}).
In the relations (\ref{D1})-(\ref{D6}) we set
\bn
\tilde{e}_{-\th^{}}=k_{\d-\a-2\b}^{-1}e_{-\a-2\b}^{}~,
\label{DY1}
\ed
where
\bn
e_{-a-2\b}^{}:=[e_{-\a-\b}^{},e_{-\b}^{}]_{q}~,
\qquad\qquad
e_{-a-\b}^{}:=[e_{-\a}^{},e_{-\b}^{}]_{q}~.
\label{DY2}
\ed
Using explicit relations for the  Cartan-Weyl
basis of $U_{q}(C_{2})\;(\simeq U_{q}(B_{2})))$ (see \cite{KT})
it is not difficult to calculate the right-hand sides of 
the relations (\ref{D2})-(\ref{D6}).
We obtain the result which is formulated as a definition of
the Drinfeldian $D_{q\eta}(C_2)$.
\begin{definition}
The Drinfeldian $D_{q\eta}(C_2)$ associated to
$C_{2}(\simeq B_{2})$ is the Hopf algebra 
generated by the quantum algebra
$U_{q}(C_{2})$ and the elements $k_{\d}^{\pm 1}$,
$\xi_{\d-\a-2\b}^{}$ with the defining relations:
\begin{eqnarray}
&&[k_{\d}^{\pm1},{\rm everything}]=0~,
\label{DY3}
\\[7pt]
&&k_{\a}^{}\xi_{\delta-\a-2\b}^{}k_{\a}^{-1}=
\xi_{\d-\a-2\b}^{}~,
\label{DY4}
\\[7pt]
&&k_{\b}^{}\xi_{\d-\a-2\b}^{}k_{\b}^{-1}=
q^{-(\a,\b)}\xi_{\d-\a-2\b}^{}~,
\label{DY5}
\\[7pt]
&&[e_{-\a}^{},\xi_{\d-\a-2\b}^{}]=-\eta q^{-\frac{1}{2}(\a,\b)}
\mbox{\small{$\left[\frac{(\a,\b)}{2}\right]$}}
k_{\d-\a-2\b}^{-1}e_{-\a-\b}^{2}~,
\label{DY6}
\\[7pt]
&&[e_{-\b}^{},\xi_{\d-\a-2\b}^{}]=0~,
\label{DY7}
\\[7pt]
&&[e_{\a}^{},\xi_{\d-\a-2\b}^{}]_{q}=\eta q^{\frac{3}{2}(\a,\b)}
[(\a,\b)]\mbox{\small{$\left[\frac{(\a,\b)}{2}\right]$}}
k_{\d-\a-2\b}^{-1}k_{\a+2\b}^{}e_{-\b}^{2}~,
\label{DY8}
\\[7pt]
&&[e_{\b}^{},[e_{\b}^{},[e_{\b}^{},
\xi_{\d-\a-2\b}^{}]_{q}]_{q}]_{q}=0~,
\label{DY9}
\\[7pt]
&&[[e_{\a}^{},\xi_{\d-\a-2\b}^{}]_{q},\xi_{\d-\a-2\b}^{}]_{q}=0~.
\label{DY10}
\end{eqnarray}
The Hopf structure of the Drinfeldian $D_{q\eta}(C_{2})$ is
defined by the formulas $\Delta_{q\eta}(x)=\Delta_{q}(x)$,
$S_{q\eta}(x)=S_{q}(x)$ $(x\in U_{q}(C_{2}))$ and
\begin{eqnarray}
&&
\Delta_{q\eta}(\xi_{\delta-\a-2\b}^{})=
\xi_{\d-\a-2\b}^{}\otimes 1+
k_{\d-\a-2\b}^{-1}\otimes\xi_{\d-\a-2\b^{}}
\nonumber
\\[3pt]
&&\quad
+\eta\big(k_{\d-\a-2\b}^{-1}\otimes k_{\d-\a-2\b}^{-1}\big)
\Big(e_{-\a-2\b}\otimes
\mbox{\small{$\frac{k_{\a+2\b}^{}-k_{\d-\a-2\b}^{}}{q-q^{-1}}$}}
-[(\a,\b)]e_{-\a-\b}\otimes k_{\a+\b}e_{-\b}\qquad\qquad
\nonumber
\\[3pt]
&&\qquad\qquad\qquad\qquad\qquad\quad
+(q-q^{-1})q^{\frac{3}{2}(\a,\b)}[(\a,\b)]
\mbox{\small{$\left[\frac{(\a,\b)}{2}\right]$}}\,
e_{-\a}^{}\otimes k_{\a}e_{-\b}^{2}\Big)~,
\label{DY11}
\end{eqnarray}
\begin{eqnarray}
&&S_{q\eta}(\xi_{\delta-\a-2\b}^{})=
-k_{\delta-\a-2\b}^{}\xi_{\d-\a-2\b}^{}
\nonumber
\\[3pt]
&&\qquad\qquad\quad
+\eta\Big(e_{-\a-2\b}^{}
\mbox{\small{$\frac{k_{\a+2\b}^{}k_{d-\a-2\b}^{-1}-1}{q-q^{-1}}$}}
-[(\a,\b)]\,e_{-\a-\b}^{}e_{-\b}^{}
\nonumber
\\[3pt]
&&\qquad\qquad\qquad
-(q-q^{-1})q^{\frac{1}{2}(\a,\b)}[(\a,\b)]\,
\mbox{\small{$\left[\frac{(\a,\b)}{2}\right]$}}
e_{-\a}^{}e_{-\b}^{2}\Big)
k_{\a+2\b}^{-1}k_{\d-\a-2\b}^{}~,\qquad\qquad
\label{DY12}
\end{eqnarray}
Here we use the standard notation
$[a]:=(q^{a}-q^{-a})/(q-q^{-1})$.
\end{definition}

At the limit $q=1$ we obtain the Yangian
$Y_{\eta}(C_{2})'=D_{q=1,\eta}(C_{2})$ with the central element
$h_{\d}$. We formulate this result as the following proposition.
\begin{proposition}
The Yangian $Y_{\eta}(C_2)$ (as an associative algebra Hopf
algebra over $\CC[[\eta]]$ is generated by the algebra
$U(C_{2})$ and the elements $h_{\d}^{}$, $\xi_{\d-\a-2\b}^{}$
with the defining relations:
\begin{eqnarray}
&&[h_{\d}^{},{\rm everything}]=0~,
\label{DY13}
\\[7pt]
&&[h_{\a}^{},\xi_{\delta-\a-2\b}^{}]=\xi_{\d-\a-2\b}^{}~,
\label{DY14}
\\[7pt]
&&[h_{\b}^{},\xi_{\d-\a-2\b}^{}]=(\a,\b)\xi_{\d-\a-2\b}^{}~,
\label{DY15}
\\[7pt]
&&[e_{-\a}^{},\xi_{\d-\a-2\b}^{}]=
-\eta\,\frac{1}{2}\,(\a,\b)\,e_{-\a-\b}^{2}~,
\label{DY16}
\\[7pt]
&&[e_{-\b}^{},\xi_{\d-\a-2\b}^{}]=0~,
\label{DY17}
\\[7pt]
&&[e_{\a}^{},\xi_{\d-\a-2\b}^{}]=
\eta\,\frac{1}{2}\,(\a,\b)^{2}e_{-\b}^{2}~,
\label{DY18}
\\[7pt]
&&[e_{\b}^{},[e_{\b}^{},[e_{\b}^{},\xi_{\d-\a-2\b}^{}]]]=0~,
\label{DY19}
\\[7pt]
&&[[e_{\a}^{},\xi_{\d-\a-2\b}^{}],\xi_{\d-\a-2\b}^{}]=0~.
\label{DY20}
\end{eqnarray}
The nontrivial coproduct $\Delta_{\eta}$ and antipode $S_{\eta}$
for the element $\xi_{\d-\a-2\b}$ is given by the formulas
\begin{eqnarray}
\Delta_{\eta}(\xi_{\delta-\a-2\b}^{})&\!\!\!=&\!\!\!
\xi_{\d-\a-2\b}^{}\otimes1+1\otimes\xi_{\d-\a-2\b}^{}
\nonumber
\\[3pt]
&&\!\!\!\!\!\!\!
+\,\eta\biggl(e_{-\a-2\b}^{}\otimes
\Big(h_{\a+2\b}^{}-\frac{h_{\d}}{2}\Big)
-(\a,\b)\,e_{-\a-\b}^{}\otimes e_{-\b}^{}\biggr)~,\qquad\qquad
\label{DY21}
\end{eqnarray}
\bn
S_{\eta}(\xi_{\delta-\a-2\b}^{})=-\xi_{\d-\a-2\b}^{}
+\eta\biggl(e_{-\a-2\b}^{}
\Big(h_{\a+2\b}^{}-\frac{h_{\d}}{2}\Big)
-(\a,\b)e_{-\a-\b}^{}e_{-\b}^{}\biggr).
\label{DY22}
\ed
\label{P1}
\end{proposition}

\vspace{5pt}
\noindent
{\large {\it 1. The Drinfeldian $D_{q\eta}(G_{2})$ and
the Yangian $Y_{\eta}(G_{2})$}}.
In the case of the Lie algebra $g=G_{2}$ there are two
positive simple roots $\a$ and $\b$, and the maximal
positive root is $\theta=2\a+3\b$.
The Dynkin diagram  of the corresponding affine Lie algebra 
$G_{2}^{(1)}$ is presented by the picture \cite{K}

{\thicklines
\vskip 15pt
\begin{picture}(500,20)
\put(208,9){\line(1,0){32}}
\put(169,6){\line(1,0){32}}
\put(165,6){\circle{8}}
\put(205,6){\circle{8}}
\put(245,6){\circle{8}}
\put(238,13){\line(1,-2){5}}
\put(238,-1){\line(1,2){5}}
\put(209,6){\line(1,0){32}}
\put(208,3){\line(1,0){32}}
\put(155,14){$\d\!-\!\th$}
\put(201,14){$\a$}
\put(242,14){$\b$}
\end{picture}}
\vskip 8pt
\centerline{\footnotesize Fig.2. Dynkin diagram of the Lie
algebra $G_{2}^{(1)}$}

\vskip 10pt
\noindent
The quantum Hopf algebra $U_q(G_{2})$ is generated by
the elements $k_{\a}^{\pm 1}$, $k_{\b}^{\pm 1}$,
$e_{\pm\a}^{}$, $e_{\pm\b}^{}$ with the defining relations
(\ref{A8})-(\ref{A13}) (see also details in \cite{KT}).
In the relations (\ref{D1})-(\ref{D6}) we set
\bn
\tilde{e}_{-\th}^{}=k_{d-2\a-3\b}^{-1}e_{-2\a-3\b}^{}~,
\label{DY23}
\ed
where
\begin{eqnarray}
&&e_{-2a-3\b}^{}:=[e_{-\a-2\b}^{},e_{-\a-\b}^{}]_{q}~,
\qquad\quad
e_{-a-\b}^{}:=[e_{-\b}^{},e_{-\a}^{}]_{q}~,
\nonumber
\\[5pt]
&&e_{-a-2\b}^{}:=[e_{-\b}^{},e_{-\a-\b}^{}]_{q}~,
\qquad\qquad\quad
e_{-a-3\b}^{}:=[e_{-\b}^{},e_{-\a-2\b}^{}]_{q}.
\label{DY24}
\end{eqnarray}
Using explicit relations for the  Cartan-Weyl
basis of $U_{q}(G_{2})$ (see \cite{KT}) we can calculate
the right-hand sides of the relations (\ref{D2})-(\ref{D6}).
We obtain the result which is formulated as a definition of
the Drinfeldian $D_{q\eta}(G_2)$.
\begin{definition}
The Drinfeldian $D_{q\eta}(G_2)$ associated to
$G_{2}(\simeq B_{2})$ is the Hopf algebra 
generated by the quantum algebra
$U_{q}(G_{2})$ and the elements $k_{\d}^{\pm 1}$,
$\xi_{\d-2\a-3\b}^{}$ with the defining relations:
\begin{eqnarray}
&&[k_{\d}^{\pm1},{\rm everything}]=0~,
\label{DY25}
\\[7pt]
&&k_{\a}^{}\xi_{\delta-2\a-3\b}^{}k_{\a}^{-1}=
q^{-(\a,\b)}\xi_{\d-2\a-3\b}^{}~,
\label{DY26}
\\[7pt]
&&k_{\b}^{}\xi_{\d-2\a-3\b}^{}k_{\b}^{-1}=
\xi_{\d-2\a-3\b}^{}~,
\label{DY27}
\\[7pt]
&&[e_{-\a}^{},\xi_{\d-2\a-3\b}^{}]=-\eta(q-q^{-1})
q^{-(\a,\b)}bc\,k_{\d-2\a-3\b}^{-1}e_{-\a-\b}^{3}~,
\label{DY28}
\\[7pt]
&&[e_{-\b}^{},\xi_{\d-2\a-3\b}^{}]=
-\eta\,q^{-\frac{1}{3}(\a,\b)}
ab^{-1}c\,k_{\d-2\a-3\b}^{-1}e_{-\a-2\b}^{2}~,
\label{DY29}
\\[7pt]
&&[e_{\b}^{},\xi_{\d-2\a-3\b}^{}]_{q}=
-\eta\,q^{\frac{1}{3}(\a,\b)}
ab\,k_{\d-2\a-3\b}^{-1}k_{\b}^{}e_{-\a-\b}^{2}~,
\label{DY30}
\\[7pt]
&&[e_{\a}^{},[e_{\a}^{},\xi_{\d-2\a-3\b}^{}]_{q}]_{q}=
\eta(q-q^{-1})\,q^{-(\a,\b)}abcd\,k_{\d-2\a-3\b}^{-1}
k_{\a}^{-2}e_{-\b}^{3}~,
\label{DY31}
\\[7pt]
&&[[e_{\a}^{},\xi_{\d-\a-2\b}^{}]_{q},\xi_{\d-\a-2\b}^{}]_{q}
=\eta^{2}\,q^{-(\a,\b)}ab^{-1}c^{2}d\,
k_{\d-2\a-3\b}^{-2}k_{\a}^{-1}e_{-\a-2\b}^{3}~.\qquad
\label{DY32}
\end{eqnarray}
The Hopf structure of the Drinfeldian $D_{q\eta}(G_{2})$ is
defined by the formulas $\Delta_{q\eta}(x)=\Delta_{q}(x)$,
$S_{q\eta}(x)=S_{q}(x)$ ($x\in U_{q}(G_{2})$) and
\begin{eqnarray}
&&
\Delta_{q\eta}(\xi_{\delta-2\a-3\b}^{})=
\xi_{\d-2\a-3\b}^{}\otimes 1+
k_{\d-2\a-3\b}^{-1}\otimes\xi_{\d-2\a-3\b}^{}
\nonumber
\\[4pt]
&&\quad
+\eta\,(k_{\d-2\a-3\b}^{-1}\otimes k_{\d-2\a-3\b}^{-1})
\biggl(e_{-2\a-3\b}^{}\otimes
\mbox{\small{$\frac{k_{2\a+3\b}^{}-k_{\d-2\a-3\b}^{}}{q-q^{-1}}$}}
\nonumber
\\[4pt]
&&\qquad
-q^{-\frac{1}{3}(\a,\b)}ae_{-\a-2\b}^{}\otimes
k_{\a+2\b}^{}e_{-\a-\b}^{}+
q^{\frac{4}{3}(\a,\b)}\Big(1-(q-q^{-1})b\Big)
e_{-\a-3\b}^{}\otimes k_{\a+3\b}^{}e_{-\a}^{})\quad
\nonumber
\\[4pt]
&&\qquad
+(q-q^{-1})a\,\Big(q^{(\a,\b)}a\,e_{-\b}^{}e_{-\a-2\b}^{}\otimes
k_{\a+2\b}^{}e_{-\a}^{}+q^{\frac{1}{3}(\a,\b)}b\,e_{-\b}^{}\otimes
k_{\b}e_{-\a-\b}^{2}\Big)\qquad\quad
\nonumber
\\[4pt]
&&\qquad
-(q-q^{-1})^{2}q^{\frac{7}{3}(\a,\b)}a^{2}b\,
e_{-\b}^{2}\otimes k_{\b}^{2}e_{-\a-\b}^{}e_{-\a}^{}
\nonumber
\\[4pt]
&&\qquad
+(q-q^{-1})^{3}q^{4(\a,\b)}a^{2}bc\,
e_{-\b}^{3}\otimes k_{\b}^{3}e_{-\a}^{2})\biggr)~,
\end{eqnarray}
\begin{eqnarray}
&&
S_{q\eta}(\xi_{\delta-2\a-3\b}^{})=
-k_{\delta-2\a-3\b}^{}\xi_{\d-2\a-3\b}^{}
\nonumber
\\[4pt]
&&\quad
+\eta\,\biggl(e_{-2\a-3\b}^{}
\mbox{\small{$\frac{k_{2\a+3\b}k_{\d-2\a-3\b}^{-1}-1}{q-q^{-1}}$}}
+q^{\frac{2}{3}(\a,\b)}a\,e_{-\a-3\b}^{}e_{-\a}^{}
+q^{-\frac{1}{3}(\a,\b)}a\,e_{-\a-2\b}^{}e_{-\a-\b}^{}
\nonumber
\\[4pt]
&&\qquad\quad
+(q-q^{-1})ab\,\Big(q^{\frac{2}{3}(\a,\b)}
e_{-\a-3\b}^{}e_{-\a}^{}-
q^{-\frac{1}{3}(\a,\b)}e_{-\b}^{}e_{-\a-\b}^{2}\Big)
\nonumber
\\[4pt]
\cr
&&\qquad\quad
+(q-q^{-1})^{2}a^{2}b\,
\Big(q^{\frac{1}{3}(\a,\b)}e_{-\b}^{}e_{-\a-2\b}^{}e_{-\a}^{}
-q^{-\frac{1}{3}(\a,\b)}e_{-\b}^{2}e_{-\a-\b}^{}e_{-\a}^{}\Big)
\nonumber
\\[4pt]
\cr
&&\qquad\quad
-(q-q^{-1})^{3}a^{2}bc\,e_{-\b}^{3}e_{-\a}^{2}\biggr)
k_{2\a+3\b}^{-1}k_{\d-2\a-3\b}^{}~.
\label{DY34}
\end{eqnarray}
\end{definition}
Here we use the notations
\begin{eqnarray}
&&a:=\frac{q^{(\a,\b)}-q^{-(\a,\b)}}{q-q^{-1}}~,
\qquad\qquad
d:=\frac{q^{2(\a,\b)}-q^{-2(\a,\b)}}{q-q^{-1}}~,
\nonumber
\\[7pt]
&&b:=\frac{q^{\frac{2}{3}(\a,\b)}-
q^{-\frac{2}{3}(\a,\b)}}{q-q^{-1}}~,\qquad\qquad
c:=\frac{q^{\frac{1}{3}(\a,\b)}-
q^{-\frac{1}{3}(\a,\b)}}{q-q^{-1}}~,
\label{DY35}
\end{eqnarray}

At the limit $q=1$ we obtain the Yangian
$Y_{\eta}(G_{2})'=D_{q=1,\eta}(G_{2})$ with the central element
$h_{\d}$. We formulate this result as the following proposition.
\begin{proposition}
The Yangian $Y_{\eta}(G_2)$ (as an associative algebra Hopf
algebra over $\CC[[\hbar,\eta]]$ is generated by the algebra
$U(C_{2})$ and the elements $h_{\d}^{}$, $\xi_{\d-\a-2\b}^{}$
with the defining relations:
\begin{eqnarray}
&&[h_{\d},{\rm everything}]=0~,
\label{DY36}
\\[7pt]
&&[h_{\a}^{},\xi_{\d-2\a-3\b}^{}]=(\a,\b)\xi_{\d-2\a-3\b}^{}~,
\label{DY37}
\\[7pt]
&&[h_{\b}^{},\xi_{\d-2\a-3\b}^{}]=\xi_{\d-2\a-3\b}^{}~,
\label{DY38}
\\[7pt]
&&[e_{-\a}^{},\xi_{\d-2\a-3\b}^{}]=0~,
\label{DY39}
\\[7pt]
&&[e_{-\b}^{},\xi_{\d-2\a-3\b}^{}]=
-\eta\,\frac{1}{2}\,(\a,\b)\,e_{-\a-2\b}^{2}~,
\label{DY40}
\\[7pt]
&&[e_{\b}^{},\xi_{\d-2\a-3\b}^{}]=
\eta\,\frac{2}{3}\,(\a,\b)^{2}e_{-\a-\b}^{2}~,
\label{DY41}
\\[7pt]
&&[e_{\a}^{},[e_{\a}^{},\xi_{\d-\a-2\b}^{}]]=0~,
\label{DY42}
\\[7pt]
&&[[e_{\a}^{},\xi_{\d-2\a-3\b}^{}],\xi_{\d-2\a-3\b}^{}]=
\eta^{2}\frac{1}{3}\,(\a,\b)^{3}\,e_{-\a-2\b}^{3}~,
\label{DY43}
\end{eqnarray}
The nontrivial coproduct $\Delta_{\eta}$ and antipode $S_{\eta}$
for the element $\xi_{\d-2\a-3\b}^{}$ is given by the formulas
\begin{eqnarray}
&&\Delta_{\eta}(\xi_{\delta-2\a-3\b}^{})=
\xi_{\d-2\a-3\b}^{}\otimes1+1\otimes\xi_{\d-2\a-3\b}^{}
\label{DY44}
\\[4pt]
&&\quad
+\,\eta\,\biggl(e_{-\a-2\b}^{}\otimes
\Big(h_{2\a+3\b}^{}-\frac{h_{\d}}{2}\Big)
+\,(\a,\b)\big(e_{-\a-3\b}^{}\otimes e_{-\a}^{}
-e_{-\a-2\b}^{}\otimes e_{-\a-\b}^{}\big)\biggr)\qquad
\nonumber
\end{eqnarray}
\begin{eqnarray}
&&S_{\eta}(\xi_{\delta-2\a-3\b}^{})=-\xi_{\d-2\a-3\b}^{}
\label{DY45}
\\[4pt]
&&\quad
+\eta\,\biggl(e_{-2\a-3\b}^{}\Big(h_{2\a+3\b}^{}-
\frac{h_{\d}}{2}\Big)+(\a,\b)(e_{-\a-3\b}^{}e_{-\a}^{}-
e_{-\a-2\b}^{}e_{-\a-\b}^{})\biggr)~.
\nonumber
\end{eqnarray}
\label{P2}
\end{proposition}

\setcounter{equation}{0}
\section{Super Drinfeldian and super Yangian}

It is obvious that results of the Sections 3 are easy
extended to the supercase.
In the supercase, i.e. when $g$ is a simple finite-dimensional
contragredient Lie superalgebra all the commutators and
the q-commutators of the relations (\ref{D2})-(\ref{D4})
are replaced by the supercommutators and
the q-supercommutators.
For example, the q-commutator (\ref{A10}) is replace by
the q-supercommutator
\bn
[e_{\b}^{},e_{\g}^{}]_{q}=e_{\b}e_{\g}-
(-1)^{\vartheta(\b)\vartheta(\g)}q^{(\b,\g)}e_{\g}^{}e_{\b}^{}~,
\label{SD1}
\ed
where $\vartheta(\cdot)$ is a standard parity function
($\vartheta(\g)=0$ for any even root $\g$, and
$\vartheta(\g)=1$ for any odd root $\g$).
Moreover we have to add some additional Serre relations
if they exist (see \cite{KT,Y}, for example).

Let us construct the super Drinfeldian and the super Yangian
of the superalgebra $osp(1|2)$ $(\simeq B(0,1))$ as an example
in an explicit form.
The Dynkin diagram of the affine superalgebra
$\widehat{osp}(1|2)$ ($\simeq B(0,1)^{(1)}$) is represented
by the picture

{\thicklines
\vskip 15pt
\begin{picture}(500,20)
\put(183,9){\line(1,0){32}}
\put(184,7){\line(1,0){32}}
\put(180,6){\circle{8}}
\put(220,6){\circle*{8}}
\put(213,13){\line(1,-2){5}}
\put(213,-1){\line(1,2){5}}
\put(184,5){\line(1,0){32}}
\put(183,3){\line(1,0){32}}
\put(170,15){$\d\!-\!2\a$}
\put(217,15){$\alpha$}
\end{picture}}
\vskip 8pt
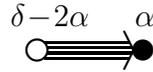
\centerline{\footnotesize Fig.3. Dynkin diagram of the Lie
superalgebra $\widehat{osp}(1|2)$ ($\simeq B(0,1)^{(1)}$).}

\vskip 10pt
\noindent
where $\a$ is the odd root ($\vartheta(\a)=1$) and
$\d-2\a$ is even one ($\vartheta(\d-2\a)=0$).

The quantum Hopf algebra $U_{q}(osp(1|2))$ is generated by
the elements $k_{\a}^{\pm 1}$, and $e_{\pm\a}^{}$ with
the defining relations
\bn
k_{\a}^{}e_{\pm\a}^{}k^{-1}_{\a}=
q^{\pm(\a,\a)}e_{\pm\a}^{}~,\qquad
[e_{\a}^{},e_{-\a}^{}]=
\frac{k_{\a}^{}-k_{\a}^{-1}}{q-q^{-1}}~,
\label{SD2}
\ed
where the brackets $[\cdot,\cdot]$ is the fermion commutator.
The Hopf structure of $U_{q}(osp(1|2))$  is given by the formulas
(\ref{A11})-(\ref{A13}) (see details in \cite{KT}).

In the relations (\ref{D1})-(\ref{D6}) we set
\bn
\tilde{e}_{-\th}^{}=e_{-\a}^{2}~.
\label{SD3}
\ed
After calculations of the right-hand sides of the relations
(\ref{D2})-(\ref{D6}) we obtain the result which is formulated
as a definition of the Drinfeldian $D_{q\eta}(osp(1|2))$.
\begin{definition}
The Drinfeldian $D_{q\eta}(osp(1|2))$ associated to
$osp(1|2)$ is the Hopf algebra generated by the quantum algebra
$U_{q}(osp(1|2))$ and the elements and $k_{\d}^{\pm1}$, 
$\xi_{\d-2\a}^{}$ with the defining relations:
\begin{eqnarray}
&&
[k_{\delta}^{\pm}, {\rm everything}]=0~,
\label{SD4}
\\[7pt]
&&
k_{\a}^{}\xi_{\d-2\a}^{}k_{\a}^{-1}=q_{\a}^{-2}\xi_{\d-2\a}^{}~,
\label{SD5}
\\[7pt]
&&
[e_{-\a}^{},\xi_{\d-2\a}^{}]=0~,
\label{SD6}
\\[7pt]
&&
[e_{\a}^{},[e_{\a}^{},[e_{\a}^{},[e_{\a}^{},[e_{\a}^{},
\xi_{\d-2\a}^{}]_q]_q]_q]_q]_q=0~,
\label{SD7}
\end{eqnarray}
\bn
\begin{array}{l}
[[e_{\a}^{},\xi_{\d-2\a}^{}]_{q},\xi_{\d-2\a}^{}]_q=
\eta^2[(\alpha,\alpha)]\left([\alpha,\alpha)]
e_{-\alpha}^4e_{\alpha}+\left[\frac{(\alpha,\alpha)}{2}\right]
[h_{\alpha}+\frac{7}{2}(\alpha,\alpha)]e_{-\alpha}^3\right)
\label{SD8}
\\[4pt]
-\eta(q-q^{-1})\left([(\alpha,\alpha)]^2
e_{-\alpha}^2\{e_{\alpha},\xi_{\delta-2\alpha}\}-
\right.\left.\left[\frac{(\alpha,\alpha)}{2}\right]
[2(\alpha,\alpha)][h_{\alpha}+\frac{5}{2}(\alpha,\alpha)]
e_{-\alpha}\xi_{\delta-2\alpha}\right)~.
\end{array}
\nonumber
\ed
The Hopf structure of the Drinfeldian $D_{q\eta}(osp(1|2))$ is
defined by the formulas $\Delta_{q\eta}(x)=\Delta_{q}(x)$,
$S_{q\eta}(x)=S_{q}(x)$ ($x\in U_{q}(osp(1|2))$) and
\begin{equation}
\begin{array}{l}
\Delta_{q\eta}(\xi_{\d-2\a}^{})=
\xi_{\delta-2\alpha}\otimes 1+
k_{\delta-2\alpha}^{-1}\otimes\xi_{\delta-2\alpha}
\nonumber
\\[4pt]
+\eta\left(e_{-\alpha}^2\otimes\frac{k_{\alpha}^2-1}{q-q^{-1}}-
\frac{k_{\delta-2\alpha}^{-1}-1}{q-q^{-1}}\otimes e_{-\alpha}^2-
q_{\alpha}^{\frac{1}{2}}\Big[\frac{(\alpha,\alpha)}{2}\Big]\,
e_{-\alpha}\otimes k_{\alpha}e_{-\alpha}\right)~,
\end{array}
\label{SD9}
\end{equation}
\begin{equation}
S_{q\eta}(\xi_{\d-2\a}^{})=
-k_{\delta-2\alpha}\xi_{\delta-2\alpha}+
\eta q_{\alpha}\frac{q_{\alpha}^3k_{\delta}-1}{q-q^{-1}}
e_{-\alpha}^2k_{\alpha}^{-2}~,
\label{SD10}
\end{equation}
\begin{equation}
\varepsilon_{q\eta}(\xi_{\d-2\a})=0~,\qquad
\varepsilon_{q\eta}(k_{\d}^{})=\varepsilon_{q\eta}(1)=1~.
\label{SD11}
\end{equation}
Here $q_{\a}:=q^{(\a,\a)}$, $\{x,y\}:=xy+yx$,
and $[a]:=(q^{a}-q^{-a})/(q-q^{-1})$.
\end{definition}

At the limit $q=1$ we obtain the super Yangian
$Y'_{\eta}(osp(1|2))=D_{q=1,\eta}(osp(1|2))$ with
the defining relations:
\begin{equation}
[h_{\d},{\rm everything}]=0~,\qquad
[e_{\alpha},e_{-\alpha}]=h_{\alpha}~,
\label{SD12}
\end{equation}
\begin{equation}
[h_{\alpha},e_{\pm\alpha}]=
\pm(\alpha,\alpha)e_{\pm\alpha}~,\qquad
[h_{\alpha},\xi_{\delta-2\alpha}]=
-2(\alpha,\alpha)\xi_{\delta-2\alpha}~,
\label{SD13}
\end{equation}
\begin{equation}
[e_{-\alpha},\xi_{\delta-2\alpha}]=0~,\qquad
[e_{\alpha},[e_{\alpha},[e_{\alpha},[e_{\alpha},
[e_{\alpha},\xi_{\delta-2\alpha}]]]]]=0~,
\label{SD14}
\end{equation}
\begin{equation}
[[e_{\alpha},\xi_{\delta-2\alpha}],\xi_{\delta-2\alpha}]=
\eta^2\frac{(\alpha,\alpha)^2}{2}
\left(2e_{-\alpha}^4e_{\alpha}+\left(h_{\alpha}+
\frac{7}{2}(\alpha,\alpha)\right)e_{-\alpha}^3\right)
\label{SD15}
\end{equation}
with the non-trivial comultiplication $\Delta_{\eta}$ and
the antipode $S_{\eta}$ for the affine root vector
$\xi_{\delta-2\alpha}$ given by
\begin{equation}
\begin{array}{ll}
\displaystyle
\Delta_{\eta}(\xi_{\delta-2\alpha})=&
\displaystyle
\xi_{\delta-2\alpha}\otimes 1+1\otimes\xi_{\delta-2\alpha}
\nonumber
\\
&
\displaystyle
+ \ \eta\left(e_{-\alpha}^2\otimes h_{\alpha}+
\frac{1}{2}h_{\delta-2\alpha}\otimes e_{-\alpha}^2-
\frac{(\alpha,\alpha)}{2}e_{-\alpha}\otimes e_{-\alpha}\right)~,
\end{array}
\label{SD16}
\end{equation}
\begin{equation}
S_{\eta}(\xi_{\delta-2\alpha})=-\xi_{\delta-2\alpha}+
{\eta\over 2}
\left( h_{\delta}+3(\alpha,\alpha)\right)
\,e_{-\alpha}^{2} \, .
\label{SD17}
\end{equation}
The Proposition \ref{PD1} is also valid for
the super Drinfeldians and the super Yangians.

\section*{Acknowledgments}

The author is thankful to Arnold Sommerfeld Institute
for Mathematical Physics, Technical University of Clausthal,
and the Organizing Committee of the International Workshop
"Lie Theory and Its Applications in Physics II",
H.-D. Doebner, V.K. Dobrev and J. Hilgert,
for the support of his visit on this Workshop.
This work was supported by the Russian Foundation
for Fundamental Research, grant No. 96-01-01421.

\end{document}